\documentclass[11pt]{amsart}
\usepackage{amssymb,amsfonts}

\newcommand{\HH}{{\mathbb H}}
\newcommand{\CC}{{\mathbb C}}
\newcommand{\RR}{{\mathbb R}}
\newcommand{\ZZ}{{\mathbb Z}}

\renewcommand{\phi}{\varphi}

\newcommand{\spec}{\rm{spec}}
\newcommand{\Spin}{\rm{Spin}}

\newcommand{\SU}{{\rm{SU}}}

\newtheorem{thm}{Theorem}
\newtheorem{lemma}[thm]{Lemma}
\newtheorem{prop}[thm]{Proposition}

\newtheorem{remark}[thm]{Remark}
\newtheorem{remarks}[thm]{Remarks}
\newtheorem{definition}[thm]{Definition}
\newtheorem{notation}[thm]{Notation}
\newtheorem{example}[thm]{Example}
\newtheorem{conjecture}[thm]{Conjecture}

\begin{document}

\title
[Prescribing eigenvalues of the Dirac operator]
{Prescribing eigenvalues of the Dirac operator}

\author{Mattias Dahl}

\address{
Institutionen f\"or Matematik\\
Kungl Tekniska H\"ogskolan\\
100 44 Stockholm\\
Sweden
}
 
\email{dahl@math.kth.se}



\date{\today}

\begin{abstract}
In this note we show that every compact spin manifold of dimension 
$\geq 3$ can be given a Riemannian metric for which 
a finite part of the spectrum of the Dirac operator consists of 
arbitrarily prescribed eigenvalues with multiplicity 1.
\end{abstract}

\maketitle

%
%%%%%%%%%%%%%%%%%%%%%%%%%%%%%%%%%%%%%%%%%%%%%%%%%%%%%%%%%%%%%%%%%%%%%%%%%
\section{Introduction and statement of results} 
%%%%%%%%%%%%%%%%%%%%%%%%%%%%%%%%%%%%%%%%%%%%%%%%%%%%%%%%%%%%%%%%%%%%%%%%%

The Dirac operator $D$ is a formally self-adjoint first order elliptic 
differential operator acting on sections of the spinor bundle over a spin 
manifold. Standard elliptic theory tells us that if the manifold is compact 
and without boundary then the spectrum of the Dirac operator is a discrete 
subset of the real line. The relationship between the spectrum of $D$ and 
topological and geometrical invariants of the manifold is a deep and 
interesting subject. Results of different types are known.

First and foremost there is the Index Theorem of Atiyah and 
Singer which gives a way of computing the Fredholm index of the Dirac 
operator by purely topological means. Another type of results are 
statements about the Dirac spectrum for generic Riemannian metrics. 
The Index Theorem gives a topological lower bound on the multiplicity 
of the zero eigenvalue of the Dirac operator, it is conjectured that 
this lower bound is sharp for generic Riemannian metrics. 

\begin{conjecture} \label{conjecture_generic_metrics} 
On any compact spin manifold for a generic Riemannian metric the space of
harmonic spinors is not larger than it is forced to be by the index theorem.
\end{conjecture}

Conjecture \ref{conjecture_generic_metrics} is known to be true for all 
manifolds of dimension $\leq 4$ \cite{maier97} and for dimension $\geq 5$ 
for a large class of manifolds including all simply connected manifolds 
\cite{baer_dahl02}. Also it is known that for a generic metric on a 
three-dimensional manifold it holds that all eigenvalues have multiplicity 
one \cite{dahl03}. It seems reasonable to believe that this holds 
true in any dimension.

In the present paper we will be concerned with the problem of finding
specific metrics on a given spin manifold with ``non-generic'' 
properties of the Dirac spectrum. As a complement to Conjecture 
\ref{conjecture_generic_metrics} there is the following.

\begin{conjecture} \label{conjecture_harmonic_spinors} 
On any compact spin manifold there is a Riemannian metric for which 
the kernel of the Dirac operator is non-trivial. 
\end{conjecture}

From the work of Hitchin \cite{hitchin74} and B\"ar \cite{baer96} it is 
known that Conjecture \ref{conjecture_harmonic_spinors} is true for 
spin manifolds of dimension $\equiv 0,1,3,7 \pmod{8}$, although exact 
information about the dimension of the kernel is not provided by their 
methods. In \cite{seeger00} Seeger constructs explicit metrics with 
non-trivial Dirac kernel on spheres of dimension $0 \pmod{4}$.

These results can be formulated as saying that it is possible to prescribe
zero to be an eigenvalue of the Dirac operator (with some unknown 
multiplicity) and then find a Riemannian metric for which this holds. We
will show that it is possible to prescribe a finite part of the spectrum 
of the Dirac operator as eigenvalues with simple multiplicity, possibly 
with the exception of the zero eigenvalue. If Conjecture 
\ref{conjecture_generic_metrics} is true we need not make this 
exception.

Since the spectrum of the Dirac operator is always symmetric about zero
in dimensions $\not\equiv 3 \pmod{4}$ we get a slightly different formulation 
of our main result depending on the dimension of the manifold. What we will 
prove is the following.
\begin{thm} \label{main_theorem} 
Let $M$ be a compact spin manifold of dimension $n \geq 3$ and let 
$L > 0$ be a real number.
\begin{itemize}
\item
Suppose that $n \equiv 3 \pmod{4}$ and let $l_1, l_2, \dots, l_m$ be 
non-zero real numbers such that  $-L < l_1 < l_2 < \dots < l_m < L$. 
Then there is a Riemannian metric $g$ on $M$ such that
$\spec\left( D_g \right) \cap \left((-L,L) \setminus \{ 0 \} \right)$ 
consists precisely of $l_1, l_2, \dots, l_m$ as simple eigenvalues.
\item
Suppose that $n \not\equiv 3 \pmod{4}$ and let $l_1, l_2, \dots, l_m$ 
be real numbers such that $0 < l_1 < l_2 < \dots < l_m < L$. Then there 
is a Riemannian metric $g$ on $M$ such that 
$\spec\left( D_g \right) \cap \left((-L,L) \setminus \{ 0 \} \right)$ 
consists precisely of $\pm l_1, \pm l_2, \dots, \pm l_m$ as simple
eigenvalues.
\end{itemize}
\end{thm}
The idea of the proof is similar to the way Conjecture 
\ref{conjecture_harmonic_spinors} is proved in \cite{baer96}. Using the 
surgery result of \cite{baer_dahl02} we begin by constructing metrics on 
spheres for which the spectrum $D$ in a given interval consists of one 
simple eigenvalue. We can then start with any Riemannian spin manifold, 
rescale the metric and take a connected sum with finitely many spheres. 
The surgery theorem of \cite{baer_dahl02} ensures that the resulting 
manifold, which is diffeomorphic to the original manifold, has a metric 
for which the Dirac eigenvalues in a given interval are precisely the 
simple eigenvalues on the spheres. 

In \cite{colindeverdiere87} Colin de Verdier shows how to find metrics
on compact manifolds for which a finite part of the spectrum of the Laplace 
operator acting on functions is arbitrarily prescribed. 
Lohkamp \cite{lohkamp96} refines this result by showing how to
prescribe simultaneously a finite part of the Laplace spectrum, 
the volume, and certain curvature invariants. 

The author wishes to express his gratitude to Christian B\"ar for 
helpful and insightful comments.

%%%%%%%%%%%%%%%%%%%%%%%%%%%%%%%%%%%%%%%%%%%%%%%%%%%%%%%%%%%%%%%%%%%%%%%%%
\section{Preliminaries}
%%%%%%%%%%%%%%%%%%%%%%%%%%%%%%%%%%%%%%%%%%%%%%%%%%%%%%%%%%%%%%%%%%%%%%%%%

A spin manifold is always assumed to be equipped with an orientation and 
a spin structure. If $M$ is a spin manifold we denote the spinor bundle 
over $M$ by $\Sigma M$.

There are two natural endomorphisms of the spinor representation, and thus
of the fibre $\Sigma_x M$. Let $n = \dim M$. If $n$ is even we consider 
Clifford multiplication by the volume form 
$\omega = (-1)^{n/2} e_1 \cdot e_2 \cdot \dots \cdot e_{n}$.
Multiplication by $\omega$ anti-commutes with multiplication by a single
tangent vector. 

Depending on the dimension $n$ the complex spinor representation has 
either a real or a quaternionic structure, see for instance
\cite[sec.~1.7]{friedrich00}. A real structure is an anti-linear 
endomorphism $\alpha$ with $\alpha^2 = \rm{Id}$, a 
quaternionic structure is an anti-linear endomorphism $\alpha$ with 
$\alpha^2 = -\rm{Id}$. If $n \equiv 0,1,6,7 \pmod{8}$ there 
is a real structure $\alpha$ on $\Sigma_x M$, if 
$n \equiv 2,3,4,5 \pmod{8}$ there is a quaternionic structure $\alpha$ 
on $\Sigma_x M$. The endomorphism $\alpha$ commutes with Clifford 
multiplication by tangent vectors if $n \equiv 2,3,6,7 \pmod{8}$ and 
anti-commutes otherwise. Both $\omega$ and $\alpha$ extend to 
parallel fields on $M$.

\begin{prop}
Suppose $(M,g)$ is a compact Riemannian spin manifold of dimension
$n \equiv 2,3,4 \pmod{8}$. Then the eigenspaces of the Dirac operator 
are quaternionic vector spaces and thus have even complex dimension.
\end{prop}

\begin{proof} 
If $n \equiv 2,3 \pmod{8}$ there is a quaternionic structure $\alpha$ 
which commutes with Clifford multiplication by tangent vectors. Since 
$\alpha$ is parallel it commutes with the Dirac operator so eigenspaces 
of $D$ are quaternionic vector spaces. 

If $n \equiv 4 \pmod{8}$ then the composition $\omega \cdot \alpha$ is a 
parallel quaternionic structure which commutes with $D$ since $\omega$ and 
$\alpha$ are both parallel and both anti-commute with Clifford 
multiplication. The eigenspaces of $D$ are quaternionic vector spaces 
with respect to this quaternionic structure.
\end{proof}

When we talk of dimensions of eigenspaces of the Dirac operator we will 
always mean quaternionic dimension if $n \equiv 2,3,4 \pmod{8}$ and 
complex dimension otherwise. We will call an eigenvalue of $D$ 
``simple'' if it is an eigenvalue for which the corresponding eigenspace 
is one-dimensional.

\begin{prop}
Suppose $(M,g)$ is a compact Riemannian spin manifold of dimension
$n \not\equiv 3,7 \pmod{8}$. Then the spectrum of the Dirac operator
is symmetric about zero.
\end{prop}

\begin{proof}
Let $\phi$ be an eigenspinor of $D$ with eigenvalue $\lambda$.

First suppose $n$ is even. Then $\omega \cdot \phi$ is an eigenspinor 
with eigenvalue $-\lambda$ since $\omega$ anti-commutes with $D$.

Next suppose $n \equiv 1,5 \pmod{8}$. Then the real/quaternionic  
structure $\alpha$ anti-commutes with $D$ so $\alpha(\phi)$ is 
an eigenspinor with eigenvalue $-\lambda$.
\end{proof}

The Index Theorem of Atiyah and Singer relates the dimension of the 
kernel of the Dirac operator to the topology of the manifold. 

If $n$ is even the spinor bundle splits as 
$\Sigma M = \Sigma^+ M \oplus \Sigma^- M$ where $\Sigma^{\pm} M$ are
the $\pm 1$-eigenbundles for multiplication by $\omega$. The Dirac 
operator is a sum $D = D^+ \oplus D^-$ where $D^{\pm}$ maps sections
of $\Sigma^{\pm} M$ to sections of $\Sigma^{\mp} M$. The Atiyah-Singer
Index Theorem states that
$\dim \ker D^+ - \dim \ker D^- = \widehat{A}(M)$ is a topological 
invariant, the $ \widehat{A}$-genus of $M$. It follows that
$\dim \ker D = \dim \ker D^+ + \dim \ker D^- \geq |\widehat{A}(M)|$.
If $n \equiv 2 \pmod{4}$ then $\widehat{A}(M) = 0$ and the only 
conclusion to be made is that $\ker D$ is even-dimensional.

In dimensions $\equiv 1,2 \pmod{8}$ the Index Theorem tells us that
$\dim \ker D \equiv \alpha(M) \pmod{2}$ and 
$\dim \ker D^+ \equiv \alpha(M) \pmod{2}$ respectively, where 
$\alpha(M) \in \ZZ/2\ZZ$ is the topological $\alpha$-genus. It 
follows that $\dim \ker D \geq |\alpha(M)|$ in dimensions 
$\equiv 1 \pmod{8}$ and $\dim \ker D \geq 2|\alpha(M)|$ in dimensions 
$\equiv 2 \pmod{8}$.

To measure closeness of the spectra of Dirac operators on different 
Riemannian manifolds we use the following definition. Let 
$\Lambda, \epsilon > 0$. Two operators with discrete real 
spectrum are called $(\Lambda,\epsilon)$-spectral close if
\begin{itemize}
\item $\pm \Lambda$ are not eigenvalues of either operator,
\item both operators have the same total number $m$ of eigenvalues in 
the interval $(-\Lambda,\Lambda)$, and
\item if the eigenvalues in $(-\Lambda,\Lambda)$ of the two operators
are denoted by $\lambda_1 \leq \dots \leq \lambda_m$ and 
$\mu_1 \leq \dots \leq \mu_m$ respectively (each eigenvalue being 
repeated according to its multiplicity), then 
$|\lambda_i - \mu_i| < \epsilon$ for $i=1,\dots,m$.
\end{itemize}

We will consider the space of Riemannian metrics on a manifold
equipped with the $C^2$-topology, when we speak of a continuous 
family of metrics it is continuous with respect to this topology.
All arguments in the proof of Theorem 1.2 of \cite{baer_dahl02}
depend continuously on the Riemannian metric in the $C^2$-topology, 
and this gives the following generalization of that result.

\begin{thm} \label{surgerythm}
Let $M$ be a closed spin manifold equipped with a spin structure and a
continuous family of Riemannian metrics $g_x$ for $x$ in some compact 
space $X$.
Let $N\subset M$ be an embedded sphere of codimension $\geq 3$ and with 
a trivialized tubular neighborhood. Let $\tilde{M}$ be the manifold with 
spin structure obtained from $M$ by surgery along $N$.

Let $\epsilon > 0$ and $\Lambda > 0$ such that 
$\pm \Lambda \notin \spec \left( D_{g_x} \right)$ for all $x$. Then there 
exists a family of Riemannian metrics $\tilde{g}_x$, $x \in X$, on 
$\tilde{M}$ such that $D_{g_x}$ and $D_{\tilde{g}_x}$ are 
$(\Lambda,\epsilon)$-spectral close for all $x \in X$.
\end{thm}

%%%%%%%%%%%%%%%%%%%%%%%%%%%%%%%%%%%%%%%%%%%%%%%%%%%%%%%%%%%%%%%%%%%%%%%%%
\section{Prescribing eigenvalues on spheres}
%%%%%%%%%%%%%%%%%%%%%%%%%%%%%%%%%%%%%%%%%%%%%%%%%%%%%%%%%%%%%%%%%%%%%%%%%

In this section we will show how to construct Riemannian metrics on 
spheres with one Dirac eigenvalue prescribed in a given interval. The 
construction will use the surgery result of the previous section and
some special Riemannian manifolds.

\begin{prop} \label{S^n} 
Let $n \geq 3$ be an integer and let $l, L$ be real numbers with $L > 0$, 
$l \neq 0$, $l \in (-L,L)$. 
\begin{itemize}
\item
If $n \equiv 3,7 \pmod{8}$ there is a metric $g^{l,L}$ on $S^n$ for which
$\spec\left( D_{g^{l,L}} \right) \cap (-L,L)$ consists only of $l$ as 
a simple eigenvalue.
\item
If $n \not\equiv 3,7 \pmod{8}$ there is a metric $g^{l,L}$ on $S^n$ for 
which $\spec\left( D_{g^{l,L}} \right) \cap (-L,L)$ consists only of 
$\pm l$ as simple eigenvalues. 
\end{itemize}
\end{prop}

To prove the proposition we first need to construct certain Riemannian 
manifolds with harmonic spinors. They are $S^3$ with the Berger metric
and manifolds with special holonomy. In \cite[sec.~3.1]{hitchin74} 
Hitchin computes explicitly the Dirac spectrum of the family 
$g^{\rm{B}}_t$ of Berger metrics on the three-dimensional 
sphere $S^3$. In particular it is shown that for a certain parameter 
interval all eigenvalues except one are uniformly bounded away from 
zero, and the remaining simple eigenvalue changes sign when the parameter 
runs through the interval. For the existence of compact manifolds with 
special holonomy as required below we refer to \cite{joyce00}.

We define the following manifolds to use as building blocks:
\begin{itemize}
\item $(V^1,g^1)$, the 1-dimensional circle with the non-bounding spin 
structure. This has $\dim_{\CC} \ker D_{g^1} = 1$. 
\item $(V^3,g^3)$, where $V^3$ is the 3-sphere and 
$g^3 = g^{\rm{B}}_{t_0} $ is the Berger 
metric for the parameter value $t_0$ for which $\dim_{\HH} \ker D_{g^3} = 1$.
\item $(V^4,g^4)$, a compact 4-manifold with holonomy $\SU(2)$. This has
$\dim_{\HH} \ker D_{g^4} = 1$, see \cite[p.~59]{wang89}, 
\cite[Thm.~3.6.1]{joyce00}.
\item $(V^6,g^6)$, a compact 6-manifold with holonomy $\SU(3)$. Then 
$\dim_{\CC} \ker D_{g^6} = 2$.
\item $(V^7,g^7)$, a compact 7-manifold with holonomy 
$\rm{G}_2$. In this case we have $\dim_{\CC} \ker D_{g^7} = 1$.
\item $(V^8,g^8)$, a compact 8-dimensional ``Bott-manifold'' with 
holonomy $\Spin (7)$. This has $\dim_{\CC} \ker D_{g^8} = 1$.
\item $(V^{10},g^{10})$, a compact 10-manifold with holonomy $\SU(5)$. 
Then $\dim_{\HH} \ker D_{g^{10}} = 1$.
\end{itemize}

The manifold $(V^8,g^8)$ plays a special role, we are going to multiply
with this manifold to increase dimensions. From 
\cite[Rem.~4, p.~12]{hitchin74} it follows that
\begin{equation} \label{add8}
\dim \ker D_{g + g^8} = \dim \ker D_{g}
\end{equation}
for any manifold $(V,g)$, where $g + g^8$ is the product metric on 
$V \times V^8$.

\begin{lemma} \label{boundarywithharmonic} Suppose $n \geq 3$.
\begin{itemize}
\item
For $n \equiv 3,7 \pmod{8}$ there is a manifold $(W^n,h^n)$ such 
that $W^n$ is a spin boundary and $\dim \ker D_{h^n} = 1$.
\item
For $n \not\equiv 3,7\pmod{8}$ there is a manifold $(W^n,h^n)$ such 
that $W^n$ is a spin boundary and $\dim \ker D_{h^n} = 2$. 
\end{itemize}
\end{lemma}

The property that $W^n$ is a spin boundary is equivalent to the sphere 
$S^n$ being spin bordant to $W^n$. In the following we use the formulas 
in \cite[Rem.~4, p.~12]{hitchin74} to find the dimension of the kernel 
of the Dirac operator on a product manifold.
 
\begin{proof}[Proof of Lemma \ref{boundarywithharmonic}.]
If $n = 3 + 8p$ we set $W^n := V^3 \times (V^8)^p$ with $h^n$ the 
product metric. Since $V^3 = S^3$ is a spin boundary we have that $W^n$ is 
a spin boundary and from (\ref{add8}) it follows that 
$\dim \ker D_{h^n} = 1$.

If $n = 7 + 8p$ we set $W^n := V^7 \times (V^8)^p$ with the product 
metric. Any compact spin manifold of dimension 7 is a spin boundary 
\cite[p.~92]{lawson_michelsohn89}, so $W^n$ is also a spin boundary.
From (\ref{add8}) it follows that the product metric has 
$\dim \ker D_{h^n} = 1$.

If $n=8p$ we take $W^n$ as the disjoint union 
$(V^8)^p + \overline{(V^8)^p}$, where the overline means taking 
the opposite orientation. Then $W^n$ is the boundary of the cylinder 
$(V^8)^p \times [0,1]$, from (\ref{add8}) we have $\dim \ker D = 1$ on
each component of $W^n$, so $\dim \ker D_{h^n} = 2$. 
 
If $n = 1 + 8p$ we set 
$W^n := V^1 \times (V^8)^p + \overline{V^1 \times (V^8)^p}$,
if $n = 2 + 8p$  we set 
$W^n := V^{10} \times (V^8)^{p-1} + \overline{
V^{10} \times (V^8)^{p-1}}$,
and if $n = 4 + 8p$ we set 
$W^n := V^4 \times (V^8)^p + \overline{V^4 \times (V^8)^p}$. 
In these cases $W^n$ has the required properties for the same reasons 
as when $n=8p$. 

If $n = 5 + 8p$ we set $W^n := V^1 \times V^4 \times (V^8)^p$. 
Since $V^4$ has quaternionic eigenspaces and $W^n$ has not we get
$\dim \ker D_{h^n} = 2$ for the product metric. If $n = 6 + 8p$ we put 
$W^n := V^6 \times (V^8)^p$ equipped with the product metric. The product 
metric has $\dim \ker D_{h^n} = 2$. In dimensions $n=5,6$ all compact spin 
manifolds are boundaries \cite[p.~92]{lawson_michelsohn89} so the 
manifolds $W^n$ are spin boundaries as well.
\end{proof}

We now prove Proposition \ref{S^n}. The low dimensions $n=3,4$ must 
be treated separately since the general argument involves breaking down 
a bordism into elementary bordisms corresponding to surgeries of 
codimension at least three, and this requires $n \geq 5$. 

\begin{proof}[Proof of Proposition \ref{S^n}.]
First assume $n=3$. Rescale the Berger metrics $g^{\rm{B}}_t$
so that the interval $(-L,L)$ contains only one Dirac eigenvalue
for $t$ close to $t_0$. Choose the parameter $t$ near $t_0$ to 
get a single (positive or negative) eigenvalue with absolute value less 
than $|l|$. Rescaling 
again increases this eigenvalue to $l$ and we have a metric $g^{l,L}$
with the required properties.

Next assume $n=4$. The manifold $V^1 \times V^3 = S^1 \times S^3$ 
equipped with the product metric $g^1 + g^3$ has 
$\dim_{\HH} \ker D_{g^1 + g^3} = 2$, see \cite[Rem.~4, p.~12]{hitchin74}.
Rescaling this metric we can assume that the interval $(-3L,3L)$ 
contains no further eigenvalues. Performing surgery on the circle 
$S := S^1 \times \{\rm{pt.}\} \subset S^1 \times S^3$ gives the 
sphere $S^4$. Since $S$ has codimension 3 the surgery result in 
Theorem \ref{surgerythm} (applied with $\Lambda$ large and $\epsilon$ 
small) tells us that there is a metric on $S^4$ for which the Dirac 
operator has eigenspaces of total dimension $2$ with eigenvalues 
between $-l$ and $l$ and no further eigenvalues in the interval 
$(-2L,2L)$. Since the Dirac spectrum is symmetric in dimension 4 
the small eigenvalues must be either $\pm l'$, $|l'| \leq |l|$, 
with multiplicity $1$, or $0$ with multiplicity 2. If we happen 
to get $0$ as an eigenvalue we can according to 
\cite[Thm.~1.3]{maier97} make a $C^2$-small perturbation of the 
metric so as to get an arbitrarily small non-zero pair $\pm l'$ as 
the first eigenvalues of the Dirac operator while the remaining 
eigenvalues are outside $(-L,L)$. (Note that the Dirac eigenvalues
depend continuously on the metric in the $C^1$-topology
\cite[Prop.~7.1]{baer96}.) 
By rescaling this metric we get a metric $g^{l,L}$ on $S^4$ with 
the required properties.

For the rest of the proof we assume $n \geq 5$.
 
Suppose $n \equiv 3,7 \pmod{8}$. From Lemma \ref{boundarywithharmonic} 
we know that $S^n$ is spin bordant to a manifold $(W^n,h^n)$ with 
$\dim \ker D_{h^n} = 1$. By rescaling $h^n$ we can assume that the 
interval $(-3L,3L)$ contains no non-zero eigenvalues of $D_{h^n}$. 
Since the sphere is simply connected and $n \geq 5$ there is a 
sequence of surgeries of codimension at least three on $W^n$ which 
will produce the sphere, see \cite[Proof of Thm.~B]{gromov_lawson80}. 
Successive applications of Theorem \ref{surgerythm} (with $\Lambda$ 
large and $\epsilon$ small enough) gives a metric on $S^n$ with one 
simple eigenvalue close to 0 (smaller than $|l|$), and no further 
eigenvalues in $(-2L,2L)$. If we happen to get a zero eigenvalue then 
we take a $C^2$-nearby metric with no zero eigenvalue (in this case 
referring to \cite[Thm.~3.10]{baer_dahl02}). Rescaling the metric 
(and possibly changing orientation) gives us the required metric $g^{l,L}$.

Finally suppose $n \not\equiv 3,7 \pmod{8}$. The same argument as 
in the previous case gives us a metric on $S^n$ for which the Dirac 
operator has eigenspaces of total dimension $2$ with eigenvalues 
between $-l$ and $l$ and no further eigenvalues in the interval 
$(-2L,2L)$. Since the Dirac spectrum is symmetric in dimensions 
$n \not\equiv 3,7 \pmod{8}$ the small eigenvalues must be either 
$\pm l'$, $|l'| \leq |l|$, with multiplicity $1$, or $0$ with 
multiplicity 2. If we get $0$ as an eigenvalue we perturb the metric 
to have an arbitrarily small pair $\pm l'$ as the first eigenvalues 
of the Dirac operator while the remaining eigenvalues are outside 
$(-L,L)$. By rescaling this metric we get a metric $g^{l,L}$ on 
$S^n$ with the required properties.
\end{proof}

%%%%%%%%%%%%%%%%%%%%%%%%%%%%%%%%%%%%%%%%%%%%%%%%%%%%%%%%%%%%%%%%%%%%%%%%%
\section{Proof of the main theorem}
%%%%%%%%%%%%%%%%%%%%%%%%%%%%%%%%%%%%%%%%%%%%%%%%%%%%%%%%%%%%%%%%%%%%%%%%%

We can now prove Theorem \ref{main_theorem}.

{\bf Theorem~\ref{main_theorem}.}
{\em
Let $M$ be a compact spin manifold of dimension $n \geq 3$ and let 
$L > 0$ be a real number.
\begin{itemize}
\item
Suppose that $n \equiv 3,7 \pmod{8}$ and let $l_1, l_2, \dots, l_m$ be 
non-zero real numbers such that  $-L < l_1 < l_2 < \dots < l_m < L$. 
Then there is a Riemannian metric $g$ on $M$ such that
$\spec\left( D_g \right) \cap \left((-L,L) \setminus \{ 0 \} \right)$ 
consists precisely of $l_1, l_2, \dots, l_m$ as simple eigenvalues.
\item
Suppose that $n \not\equiv 3,7 \pmod{8}$ and let $l_1, l_2, \dots, l_m$ 
be real numbers such that $0 < l_1 < l_2 < \dots < l_m < L$. Then there 
is a Riemannian metric $g$ on $M$ such that 
$\spec\left( D_g \right) \cap \left((-L,L) \setminus \{ 0 \} \right)$ 
consists precisely of $\pm l_1, \pm l_2, \dots, \pm l_m$ as simple
eigenvalues.
\end{itemize}
}

\begin{proof}
Start with a metric on $M$ for which the dimension of the kernel of 
the Dirac operator is minimal among all Riemannian metrics on $M$. 
A generic choice of metric will have this property, see the discussion 
in \cite[Sec.~3]{baer_dahl02}. Rescale this metric to get a metric 
$g_0$ which has no non-zero eigenvalues in the interval $(-3L,3L)$.

For each $i=1, \dots, m$ take a metric $g^{l_i,3L}$ on $S^n$ as
provided by Proposition \ref{S^n}. Rescale these metrics as
$t_i g^{l_i,3L}$, where the $t_i$ are real-valued parameters close 
to 1. This gives metrics with one eigenvalue $\lambda_i(t_i)$ being 
a monotone function of $t_i$ passing through the value $l_i$ when 
$t_i = 1$ and no further eigenvalues in $(-2L,2L)$. 
The non-zero eigenvalues in $(-2L,2L)$ of the disjoint union 
\begin{equation} \label{disjointunion}
(M,g_0) \sqcup (S^n, t_1 g^{l_1, 3L} ) \sqcup \dots 
\sqcup (S^n, t_m g^{l_m, 3L}) 
\end{equation}
are then precisely $\lambda_1(t_1), \dots , \lambda_m(t_m)$.
Define the map 
$F: (t_1, \dots, t_m) \mapsto (\lambda_1(t_1), \dots , \lambda_m(t_m))$.
Since the functions $t_i \mapsto \lambda_i(t_i)$ are monotone for 
$t_i$ near $1$ the map $F$ restricted to a sphere $S = \partial B$ 
bounding a small ball $B \subset \RR^m$ containing $(1,\dots ,1)$ 
has degree $\pm 1$ as a map to $\RR^m \setminus (l_1, \dots, l_m)$.

Theorem \ref{surgerythm} (with sufficiently small $\epsilon$ 
and sufficiently large $\Lambda$) tells us that there is a 
family of metrics parametrized by $B$ on the connected sum
$$
M \# 
\underbrace{S^n \# \dots \# S^n}_{m \ \rm{copies} }
= M
$$ 
whose Dirac operators are $(\Lambda,\epsilon)$-spectral close to the 
Dirac operators on the disjoint union (\ref{disjointunion}) for
$(t_1, \dots, t_m) \in B $. Since the zero eigenvalue of $D_{g_0}$ 
was assumed to have minimal multiplicity the zero eigenvalue on the 
connected sum must have the same minimal multiplicity. Thus the 
summand $(M,g_0)$ cannot contribute any small non-zero eigenvalue 
to the connected sum metric.

We choose $B$ and $\epsilon$ small enough so that the non-zero 
eigenvalues $\tilde{\lambda}_i (t_1, \dots, t_m)$ are simple and 
consistently numbered for $(t_1, \dots, t_m) \in B $. We then get a 
continuous map $\tilde{F}: (t_1, \dots, t_m) \mapsto 
(\tilde{\lambda}_1, \dots, \tilde{\lambda}_m)$ which just as $F$ has 
degree $\pm 1$ as a map from $S$ to $\RR^m \setminus (l_1, \dots, l_m)$.
It follows that $\tilde{F}$ must map some point in $B$ to 
$(l_1, \dots, l_m)$, so we get a metric with the prescibed eigenvalues 
in $(-L,L) \setminus \{ 0 \}$.
\end{proof}

If Conjecture \ref{conjecture_generic_metrics} is true for a manifold
of dimension $\geq 3$ we can choose the initial metric 
to have the multiplicity of the zero eigenvalue given by the index 
theorem, and then prescribe arbitrarily a finite number of non-zero 
eigenvalues of multiplicity 1.

It seems reasonable to conjecture that one can prescribe finitely many 
Dirac eigenvalues with arbitrary finite multiplicities. In short one 
could state the following conjecture.

\begin{conjecture}
Except for the algebraic constraints on the spinor bundle 
(giving quaternionic eigenspaces and symmetric spectrum) and the 
constraints on the zero eigenvalue coming from the Atiyah-Singer index 
theorem it is possible to find a Riemannian metric on any compact spin 
manifold with a finite part of its Dirac spectrum arbitrarily prescribed.
\end{conjecture}

The reason why the method of this paper does not work for prescribing
a multiple eigenvalue is that there is then no well-defined numbering 
for the eigenvalue functions as functions on the parameter space.

%%%%%%%%%%%%%%%%%%%%%%%%%%%%%%%%%%%%%%%%%%%%%%%%%%%%%%%%%%%%%%%%%%%%%%%%%
\end{document}